\documentclass[12pt]{amsart}

\newtheorem{thm}{Theorem}[section]

\newtheorem{conj}[thm]{Conjecture}
\newtheorem{issue}{Issue}

\theoremstyle{definition}

\theoremstyle{remark}

\newcommand{\NR}{{{}^\N\R}}
\newcommand{\Iff}{\Leftrightarrow}
\newcommand{\mypar}[1]{\par\medskip\noindent\textbf{#1.}}

\renewcommand{\b}{{\mathfrak b}}

\newcommand{\p}{{\mathfrak p}}

\newcommand{\upannouncement}[1]{[\S\ref{#1} above]}

\newcommand{\cP}{\mathcal{P}}
\newcommand{\M}{\mathcal{M}}

\newcommand{\cov}{\mathsf{cov}}
\newcommand{\add}{\mathsf{add}}

\newcommand{\non}{\mathsf{non}}

\newcommand{\R}{\mathbb{R}}

\newcommand{\w}{\omega}
\renewcommand{\b}{\mathfrak{b}}

\newcommand{\bq}{\begin{quote}}
\newcommand{\eq}{\end{quote}}
\renewcommand{\O}{\mathcal{O}}
\newcommand{\B}{\mathcal{B}}
\newcommand{\BG}{\B_\Gamma}

\newcommand{\BO}{\B_\Omega}

\newcommand{\sone}{\mathsf{S}_1}    \newcommand{\sfin}{\mathsf{S}_{fin}}
\newcommand{\ufin}{\mathsf{U}_{fin}}

\newcommand{\cF}{\mathcal{F}}

\newcommand{\cU}{\mathcal{U}}

\newcommand{\naturals}{{\mathbb N}}
\newcommand{\N}{\naturals}

\newcommand{\sbst}{\subseteq}
\newcommand{\by}[2]{\par\hfill\emph{#1}, #2}

\newcommand{\Tau}{\mathrm{T}}
\newcommand{\CE}{\textsc{CE}}

\newcommand{\be}{\begin{enumerate}}
\newcommand{\ee}{\end{enumerate}}
\newcommand{\bi}{\begin{itemize}}
\newcommand{\ei}{\end{itemize}}

\newcommand{\general}{\small\vfill\par\noindent\hrulefill\par
\noindent\textbf{Previous issues.}
The first issues of this bulletin,
which contain general information (first issue),
basic definitions, research announcements, and open problems (all issues) are available online,
on \arx{math.GN/$x$}, where $x$ is \texttt{0301011}, \texttt{0302062}, \texttt{0303057},
\texttt{0304087}, \texttt{0305367}, and \texttt{0312140}, respectively, for issues number $1$ to $6$.\\[0.1cm]
\textbf{Contributions.}
Please submit your contributions (announcements, discussions, and open problems)
by e-mailing us. It is preferred to write them
in \LaTeX{}.
The authors are urged to use as standard notation as possible, or otherwise give
the definitions or a reference to where the notation is explained.
Contributions to this bulletin would not require any transfer of copyright,
and material presented here can be published elsewhere.\\[0.1cm]
\textbf{Subscription.}
To receive this bulletin (free) to your
e-mailbox, e-mail us.
}

\newcommand{\ArxPaper}[5]{\subsection{#3}{#5}\par\hfill{\arx{#1}}\par\hfill\emph{#4}\par\hfill{#2}}
\newcommand{\SPMBul}{\textbf{$\mathcal{SPM}$ Bulletin}}

\newcommand{\arx}[1]{\texttt{http://arxiv.org/abs/#1}}

\newcommand{\probmonth}{\emph{Problem of the month}}

\title[$\mathcal{SPM}$ Bulletin \textbf{7} (January 2004)]{%
$\mathcal{SPM}$ Bulletin\\[0.5cm]
Issue number 7: January 2004 \CE{}}

\begin{document}
\maketitle

\tableofcontents

\section{Editor's note}

\mypar{Open problem solved}
We are glad to announce the first solution of a \emph{Problem of the Month}
posed in the \SPMBul{}. The problem from the third issue was solved
by Lubomyr Zdomsky, a student of Taras Banakh. This solution is a part of
a large project carried by these two mathematicians, which will hopefully
be announced in this bulletin when it is finished.

\mypar{Past problems}
We also have a mew section on past problems in the \SPMBul{}.

\mypar{Paper on open problems}
We have finished writing a paper containing a significant portion of
the important problems in the field of SPM.
Email us to get a copy of the paper.

\mypar{ArXiv papers}
From now on, we will try to include announcements of
papers which are of interest to readers of this bulletin
and were recently announced in the
\emph{Mathematics ArXiv}.

We encourage all contributors to submit their papers to the Mathematics ArXiv
prior to the submission of an announcement to the \SPMBul{},
so to guarantee larger exposure for their papers. (The submission to the ArXiv
does not guarantee the inclusion of the announcement in this bulletin, so please also email us
a note concerning your submission to the ArXiv after you complete it.)
Submissions to other e-print servers can also be considered, upon request.

\medskip

Contributions to the next issue are, as always, welcome.

\medskip

\by{Boaz Tsaban}{tsaban@math.huji.ac.il}

\hfill \texttt{http://www.cs.biu.ac.il/\~{}tsaban}

\section{Research announcements}

\ArxPaper{math.LO/0311443}{miller@math.wisc.edu}
{Models in which every nonmeager set is nonmeager in a nowhere dense Cantor set}
{Maxim R.\ Burke, Arnold W.\ Miller}
{We prove that it is relatively consistent with ZFC that in any perfect Polish
space, for every nonmeager set $A$ there exists a nowhere dense Cantor set $C$ such
that $A\cap C$ is nonmeager in $C$. We also examine variants of this result
and establish a measure theoretic analog.}

\ArxPaper{math.LO/0312308}{miller@math.wisc.edu}
{The $\gamma$-Borel conjecture}
{Arnold W. Miller}
{In this paper we show that it is relatively consistent with ZFC that every
$\gamma$-set is countable while not every strong measure zero set is countable.
This answers a question of Paul Szeptycki. A set is a $\gamma$-set iff every
$\omega$-cover contains a $\gamma$-subcover. An open cover is an $\omega$-cover iff
every finite set is covered by some element of the cover. An open cover is a
$\gamma$-cover iff every element of the space is in all but finitely many elements
of the cover. Gerlits and Nagy proved that every $\gamma$-set has strong measure
zero. We also show that is consistent that every strong $\gamma$-set is countable
while there exists an uncountable $\gamma$-set. On the other hand every strong
measure zero set is countable iff every set with the Rothberger property is
countable.}
\label{BC}

\ArxPaper{math.OA/0312135}{nweaver@math.wustl.edu}
{Consistency of a counterexample to Naimark's problem}
{Charles Akemann, Nik Weaver}
{We construct a C*-algebra that has only one irreducible representation up to unitary equivalence
but is not isomorphic to the algebra of compact operators on any Hilbert space. This answers an
old question of Naimark. Our construction uses a combinatorial statement called the diamond principle,
which is known to be consistent with but not provable from the standard axioms of set theory (assuming
those axioms are consistent). We prove that the statement ``there exists a counterexample to Naimark's
problem which is generated by $\aleph_1$ elements'' is undecidable in standard set theory.}

\ArxPaper{math.LO/0312472}{}
{Comparing the automorphism group of the measure algebra with some groups
related to the infinite permutation group of the natural numbers}
{Pietro Ursino}
{We prove, by a straight construction, that the automorphism group of the
measure algebra and the subgroup of the measure preserving ones cannot be
isomorphic to the trivial automorphisms of $P(\N)/fin$.}

\ArxPaper{math.AG/0401079}{zell@math.gatech.edu}
{Quantitative study of semi-Pfaffian sets}
{Thierry Zell}
{We study the topological complexity of sets defined using Khovanskii's
Pfaffian functions, in terms of an appropriate notion of format for those sets.
We consider semi- and sub-Pfaffian sets, but more generally any definable set
in the o-minimal structure generated by the Pfaffian functions, using the
construction of that structure via Gabrielov's notion of limit sets. All the
results revolve around giving effective upper-bounds on the Betti numbers (for
the singular homology) of those sets.}

\ArxPaper{math.GN/0307225}{tsaban@math.huji.ac.il}
{$o$-bounded groups and other topological groups with strong combinatorial properties}
{Boaz Tsaban}
{We construct several topological groups with very strong combinatorial properties.
In particular, we give simple examples of subgroups of $\R$ (thus strictly $o$-bounded)
which have the Hurewicz property but are not $\sigma$-compact,
and show that the product of two 
$o$-bounded subgroups of $\NR$ may fail to be $o$-bounded,
even when they satisfy the stronger property $\sone(\BO,\BO)$.
This solves a problem of Tka\v{c}enko and Hernandez,
and extends independent solutions of Krawczyk and Michalewski
and of Banakh, Nickolas, and Sanchis.
We also construct separable metrizable groups $G$ of size continuum 
such that every countable Borel
$\w$-cover of $G$ contains a $\gamma$-cover of $G$.}

\section{Problem of the month}
Let us write $BC(P)$ for the 
Borel Conjecture for sets with property $P$, that is, 
the hypothesis that every set of reals with property $P$ is
countable.

In \upannouncement{BC}
Miller proves that $BC(\sone(\O,\O))$ implies (and is therefore equivalent to)
$BC(SMZ)$,
where SMZ stands for \emph{strong measure zero}.
The proof splits into two cases: $\aleph_1=\b$ and
$\aleph_1<\b$.
In the case $\aleph_1<\b$ Miller really shows that 
strong measure zero plus $\ufin(\O,\Gamma)$ implies $\sone(\O,\O)$.
In Theorems 14 and 19 of \cite{NSW}
it is shown that in fact, 
$$SMZ+\ufin(\O,\Gamma)\Iff\sone(\O,\O)+\ufin(\O,\Gamma)\Iff (*)_{GN},$$
where $(*)_{GN}$ is the Gerlitz-Nagy covering proeprty introduced in
\cite{GN}. 
This characterization implies that
$(*)_{GN}$ is strictly stronger than $\sone(\O,\O)$,
since it implies that 
$\non((*)_{GN})=\min\{\non(\sone(\O,\O)),\non(\ufin(\O,\Gamma))\}=
\min\{\cov(\M),\b\}(=\add(\M))$, and it is consistent that $\b<\cov(\M)$.

So we have:

\begin{thm} 
$BC((*)_{GN})\Iff BC(SMZ)$.
\end{thm}
\begin{proof}
Assume that $\aleph_1<\b$ and there exists an uncountable strong measure zero 
set $X$. As $SMZ$ is hereditary, we may assume that $|X|=\aleph_1$, and therefore
$X$ satisfies $\ufin(\O,\Gamma)$ as well, that is, $X$ satisfies $(*)_{GN}$.

Next, assumet that $\aleph_1=\b$.

Consider the collection $\Omega^{gp}$ of open $\omega$-covers $\cU$ of $X$ such that
there exists a partition $\cP$ of $\cU$ into finite sets such that
for each finite $F\sbst X$ and all but finitely many $\cF\in\cP$,
there exists $U\in\cF$ such that $F\sbst U$ \cite{coc7}.

In \cite{coc7} it is proved that $\sone(\Omega,\Omega^{gp})$
is equivalent to having $(*)_{GN}$ in all finite powers.
Now if $\aleph_1=\b$ then by \cite{ideals}
there exists an uncountable element $X$ in $\sone(\Omega,\Omega^{gp})$ 
(in particular, $X$ satisfies $(*)_{GN}$).
\end{proof}

It is a conjecture of Tomasz Weiss 
that $(*)_{GN}$ is closed under taking finite products. 
If $(*)_{GN}$ is closed under taking 
finite \emph{powers}, then $(*)_{GN}=\sone(\Omega,\Omega^{gp})$ and we
have that $BC(\sone(\Omega,\Omega^{gp}))\Iff BC(SMZ)$.
Otherwise, I do not even know whether 
$BC(\sone(\Omega,\Omega))$ implies $BC(SMZ)$
(it does if $(*)_{GN}$ implies $\sone(\O,\O)$ in all finite powers.)
In fact we need to prove the following.

\begin{conj}
If $X$ has strong measure zero and $|X|<\b$, 
then all finite powers of $X$ have strong measure zero
(equivalently, all finite powers of $X$ satisfy $\sone(\O,\O)$).
\end{conj}

This constitutes the \probmonth{}.

\by{Boaz Tsaban}{tsaban@math.huji.ac.il}

\section{Problems from earlier issues}
In this section we list the past problems posed in the \SPMBul{}.
For definitions, motivation and related results, consult the
corresponding issue.

For conciseness, we make the convention that
all spaces in question are
zero-dimentional, separable metrizble spaces.

\begin{issue}
Is $\binom{\Omega}{\Gamma}=\binom{\Omega}{\Tau}$?
\end{issue}

\begin{issue}
Is $\ufin(\Gamma,\Omega)=\sfin(\Gamma,\Omega)$?
And if not, does $\ufin(\Gamma,\Gamma)$ imply
$\sfin(\Gamma,\Omega)$?
\end{issue}

\begin{issue}
Does there exist (in ZFC) a set satisfying
$\ufin(\O,\O)$ but not $\ufin(\O,\Gamma)$?
\end{issue}
\begin{proof}[Solution]
\textbf{Yes} (Lubomyr Zdomsky, 2003).
\end{proof}

\begin{issue}
Does $\sone(\Omega,\Tau)$ imply $\ufin(\Gamma,\Gamma)$?
\end{issue}

\begin{issue}
Is $\p=\p^*$?
\end{issue}

\begin{issue}
Does there exist (in ZFC) an uncountable set satisfying $\sone(\BG,\B)$?
\end{issue}

\general


\begin{thebibliography}{00}
\bibitem{ideals}
T.\ Bartoszynski and B.\ Tsaban,
\emph{Hereditary topological diagonalizations and the Menger-Hurewicz Conjectures},
Proceedings of the American Mathematical Society,  to appear. \arx{math.LO/0208224}

\bibitem{GN}
J.\ Gerlits and Zs.\ Nagy,
\emph{Some properties of C(X), I},
Topology and its Applications \textbf{14} (1982),
151--161.

\bibitem{coc7}
Lj.\ D.\ R.\ Ko\v{c}inac and M.\ Scheepers,
\emph{Combinatorics of open covers (VII): groupability},
Fundamenta Mathematicae \textbf{179} (2003), 
131--155.

\bibitem{NSW}
A.\ Nowik, M.\ Scheepers, and T.\ Weiss,
\emph{The algebraic sum of sets of real numbers with strong measure zero sets},
J.\ Symbolic Logic \textbf{63} (1998), 301--324.

\end{thebibliography}
\end{document}